\documentclass{elsart}
\usepackage{amssymb}
\usepackage{graphicx}

\begin{document}

\begin{frontmatter}

\title{Li\'{e}nard's System and Smale's Problem}

\author{Valery A. Gaiko\thanksref{label1}}
\ead{vlrgk@yahoo.com}
\thanks[label1]{The author is grateful to the Department of Applied
Mathematical Analysis of TU~Delft (the Netherlands), the
Department of Mathematical Sciences of the University of Memphis
(USA), and the Institut des Hautes \'{E}tudes Scientifiques
(France) for their hospitality in 2003--2006.}
\address{Department~of~Mathematics
Belarusian~State~University~of~Informatics~and~Radioelectronics
L.\,Beda~Str.\,6-4,~Minsk~220040,~Belarus}

\begin{abstract}
In this paper, using geometric properties of the field rotation
parameters, we present a solution of \emph{Smale's Thirteenth
Problem} on the maximum number of limit cycles for Li\'{e}nard's
polynomial system. We also generalize the obtained result and
present a solution of \emph{Hil\-bert's Sixteenth Problem} on the
maximum number of limit cycles surrounding a singular point for an
arbitrary polynomial system. Besides, we consider a generalized
Li\'{e}nard's cubic system with three finite singularities, for
which the developed geometric approach can complete its global
qualitative analysis: in particular, it easily solves the problem
on the maximum number of limit cycles in their different
distribution. We give also an alternative proof of the main
theorem for the generalized Li\'{e}nard's system applying the
Wintner--Perko termination principle for multiple limit cycles and
discuss some other results concerning this system.
\end{abstract}

\end{frontmatter}

\section{Introduction}
\label{1}

We consider planar dynamical systems
    $$
    \dot{x}=P_{n}(x,y),
    \quad
    \dot{y}=Q_{n}(x,y),
    \eqno(1.1)
    $$
where $P_{n}(x,y)$ and $Q_{n}(x,y)$ are polynomials with real
coefficients in the real variables $x,\:y,$ and, first of all, a
special case of (1.1): classical Li\'{e}nard's polynomial system
of the form
    $$
    \dot{x}=y,
    \quad
    \dot{y}=-x+\mu_{1}\,y+\mu_{2}\,y^{2}+\mu_{3}\,y^{3}+\ldots
    +\mu_{2k}\,y^{2k}+\mu_{2k+1}\,y^{2k+1}.
    \eqno(1.2)
    $$
The main problem of qualitative theory of such systems is
\emph{Hilbert's Sixteenth Problem} on the maximum number and
relative position of their limit cycles, i.\,e., closed isolated
trajectories of (1.1). This problem was formulated as one of the
fundamental problems for mathematicians of the XX century, however
it has not been solved even in the simplest (quadratic, cubic,
etc.) cases of the polynomial systems. In this paper, we suggest a
new geometric approach to solving the problem in the case of
Li\'{e}nard's system (1.2). In this special case, it is considered
as \emph{Smale's Thirteenth Problem} becoming one of the main
problems for mathematicians of the XXI century \cite{ilyash},
\cite{smale}.
    \par
In Section~2 of this paper, applying a canonical system with field
rotation parameters and using geometric properties of the spirals
filling the interior and exterior domains of limit cycles, we
present a solution of \emph{Smale's Thirteenth Problem} for
Li\'{e}nard's polynomial system (1.2). In Section~3, by means of
the same geometric approach, we generalize the obtained result and
present a solution of \emph{Hil\-bert's sixteenth problem} on the
maximum number of limit cycles surrounding a singular point for an
arbitrary polynomial system. In Section~4, we consider a
generalized Li\'{e}nard's cubic system with three finite
singularities, for which the developed geometric approach can
complete its global qualitative analysis: in particular, it easily
solves the problem on the maximum number of limit cycles in their
different distribution. In this section, we give also an
alternative proof of the main theorem for the generalized
Li\'{e}nard's system applying the Wintner--Perko termination
principle for multiple limit cycles and discuss some other results
concerning this system.

\section{Li\'{e}nard's polynomial system}
\label{2}

System (1.2) and more general Li\'{e}nard's systems have been
studying in numerous works (see, for example, \cite{BL},
\cite{blowslloyd}, \cite{chrlloyd}, \cite{franc}, \cite{gaib},
\cite{gasullt}, \cite{kkr}--\cite{lynch}, \cite{rych}). It is easy
to see that (1.2) has the only finite singularity: an anti-saddle
at the origin. At infinity, system (1.2) for $k\geqslant1$ has two
singular points: a node at the ``ends'' of the $y$-axis and a
saddle at the ``ends'' of the $x$-axis. For studying the infinite
singularities, the methods applied in \cite{BL} for Rayleigh's and
van der Pol's equations and also Erugin's two-isocline method
developed in \cite{Gaiko} can be used. Following \cite{Gaiko}, we
will study limit cycle bifurcations of (1.2) by means of a
canonical system containing only the field rotation parameters of
(1.2). It is valid the following theorem.
    \par
    \medskip
\noindent\textbf{Theorem 2.1.}
    \emph{Li\'{e}nard's polynomial system $(1.2)$ with limit cycles
can be reduced to the canonical form}
    $$
    \dot{x}=y\equiv P,
    \quad
    \dot{y}=-x+\mu_{1}\,y+y^{2}+\mu_{3}\,y^{3}+\ldots
    +y^{2k}+\mu_{2k+1}\,y^{2k+1}\equiv Q,
    \eqno(2.1)
    $$
\emph{where $\mu_{1},$ $\mu_{3},\ldots,$ $\mu_{2k+1}$ are field
rotation parameters of $(2.1).$}
    \medskip
    \par
    \noindent\textbf{Proof.} \
Vanish all odd parameters of~(1.2),
    $$
    \dot{x}=y,
    \quad
    \dot{y}=-x+\mu_{2}\,y^{2}+\mu_{4}\,y^{4}+\ldots+\mu_{2k}\,y^{2k},
    \eqno(2.2)
    $$
and consider the corresponding equation
    $$
    \frac{dy}{dx}
    =\frac{-x+\mu_{2}\,y^{2}+\mu_{4}\,y^{4}+\ldots+\mu_{2k}\,y^{2k}}{y}
    \equiv F(x,y).
    \eqno(2.3)
    $$
Since $F(x,-y)=-F(x,y),$ the direction field of (2.3) (and the
vector field of (2.2) as well) is symmetric with respect to the
$x$-axis. It follows that for arbitrary values of the parameters
$\mu_{2},$ $\mu_{4},\ldots,$ $\mu_{2k}$ system (2.2) has a center
at the origin and cannot have a limit cycle surrounding this
point. Therefore, without loss of generality, all even parameters
of system (1.2) can be supposed to be equal, for example, to one:
$\mu_{2}=\mu_{4}=\ldots=\mu_{2k}=1$ (they could be also supposed
to be equal to zero).
    \par
To prove that the rest (odd) parameters rotate the vector field of
(2.1), let us calculate the following determinants:\\
    $$
\Delta_{\mu_{1}}=PQ'_{\mu_{1}}-QP'_{\mu_{1}}=y^{2}\geq0,
    $$
    $$
\Delta_{\mu_{3}}=PQ'_{\mu_{3}}-QP'_{\mu_{3}}=y^{2}\geq0,
    $$
    $$.\,.\,.\,.\,.\,.\,.\,.\,.\,.\,.\,.\,.\,.\,.\,.
    \,.\,.\,.\,.\,.\,.\,.\,.\,.\,.\,.\,.\,.\,.\,.\,$$
    $$
\Delta_{\mu_{2k+1}}=PQ'_{\mu_{2k+1}}-QP'_{\mu_{2k+1}}=y^{2}\geq0.
    $$
    \par
By definition of a field rotation parameter~\cite{duff}, for
increasing each of the parameters $\mu_{1},$ $\mu_{3},\ldots,$
$\mu_{2k+1},$ under the fixed others, the vector field of system
(2.1) is rotated in positive direction (counter\-clock\-wise) in
the whole phase plane; and, conversely, for decreasing each of
these parameters, the vector field of (2.1) is rotated in negative
direction (clock\-wise).
    \par
Thus, for studying limit cycle bifurcations of (1.2), it is
sufficient to consider canonical system (2.1) containing only its
odd parameters, $\mu_{1},$ $\mu_{3},\ldots,$ $\mu_{2k+1},$ which
rotate the vector field of (2.1). The theorem is proved.
    \par
By means of canonical system (2.1), let us study global limit
cycle bifurcations of (1.2) and prove the following theorem.
    \par
    \medskip
\noindent\textbf{Theorem 2.2.}
    \emph{Li\'{e}nard's polynomial system $(1.2)$ has at most
$k$~limit cycles.}
    \medskip
    \par
    \noindent\textbf{Proof.} \
According to Theorem~2.1, for the study of limit cycle
bifurcations of system (1.2), it is sufficient to consider
canonical system (2.1) containing only the field rotation
parameters of (1.2): $\mu_{1},$ $\mu_{3},\ldots,$ $\mu_{2k+1}.$
    \par
Vanish all these parameters:
    $$
    \dot{x}=y,
    \quad
    \dot{y}=-x+y^{2}+y^{4}+\ldots+y^{2k}.
    \eqno(2.4)
    $$
System (2.4) is symmetric with respect to the $x$-axis and has a
center at the origin. Let us input successively the field rotation
parameters into this system beginning with the parameters at the
highest degrees of~$y$ and alternating with their signs. So, begin
with the parameter~$\mu_{2k+1}$ and let, for definiteness,
$\mu_{2k+1}>0\!:$
    $$
    \dot{x}=y,
    \quad
    \dot{y}=-x+y^{2}+y^{4}+\ldots+y^{2k}+\mu_{2k+1}\,y^{2k+1}.
    \eqno(2.5)
    $$
In this case, the vector field of~(2.5) is rotated in positive
direction (counterclockwise) turning the origin into a nonrough
unstable focus.
    \par
Fix $\mu_{2k+1}$ and input the parameter $\mu_{2k-1}<0$ into
(2.5):
    $$
    \dot{x}=y,
    \quad
    \dot{y}=-x+y^{2}+y^{4}+\ldots+\mu_{2k-1}\,y^{2k-1}
    +y^{2k}+\mu_{2k+1}\,y^{2k+1}.
    \eqno(2.6)
    $$
Then the vector field of~(2.6) is rotated in opposite direction
(clockwise) and the focus immediately changes the character of its
stability (since its degree of nonroughness decreases and the sign
of the field rotation parameter at the lower degree of~$y$
changes) generating a stable limit cycle. Under further decreasing
$\mu_{2k-1},$ this limit cycle will expand infinitely, not
disappearing at infinity (because of the parameter $\mu_{2k+1}$ at
the higher degree of~$y).$
    \par
Denote the limit cycle by $\Gamma\!_{1},$ the domain outside the
cycle by $D_{1},$ the domain inside the cycle by $D_{2}$ and
consider logical possibilities of the appearance of other
(semi-stable) limit cycles from a ``trajectory concentration''
surrounding the origin. It is clear that, under decreasing the
parameter $\mu_{2k-1},$ a semi-stable limit cycle cannot appear in
the domain $D_{2},$ since the focus spirals filling this domain
will untwist and the distance between their coils will increase
because of the vector field rotation.
    \par
By contradiction, we can also prove that a semi-stable limit cycle
cannot appear in the domain $D_{1}.$ Suppose it appears in this
domain for some values of the parameters $\mu_{2k+1}^{*}>0$ and
$\mu_{2k-1}^{*}<0.$ Return to initial system (2.4) and change the
inputting order for the field rotation parameters. Input first the
parameter $\mu_{2k-1}<0\!:$
    $$
    \dot{x}=y,
    \quad
    \dot{y}=-x+y^{2}+y^{4}+\ldots+\mu_{2k-1}\,y^{2k-1}+y^{2k}.
    \eqno(2.7)
    $$
Fix it under $\mu_{2k-1}=\mu_{2k-1}^{*}.$ The vector field
of~(2.7) is rotated clockwise and the origin turns into a nonrough
stable focus. Inputting the parameter $\mu_{2k+1}>0$ into (2.7),
we get again system (2.6), the vector field of which is rotated
counterclockwise. Under this rotation, a stable limit cycle
$\Gamma\!_{1}$ will immediately appear from infinity, more
precisely, from a separatrix cycle of the Poincar\'{e} circle form
containing infinite singularities of the saddle and node
types~\cite{BL}. This cycle will contract, the outside spirals
winding onto the cycle will untwist and the distance between their
coils will increase under increasing $\mu_{2k+1}$ to the value
$\mu_{2k+1}^{*}.$ It follows that there are no values of
$\mu_{2k-1}^{*}<0$ and $\mu_{2k+1}^{*}>0,$ for which a semi-stable
limit cycle could appear in the domain~$D_{1}.$
    \par
This contradiction proves the uniqueness of a limit cycle
surrounding the origin in system (2.6) for any values of the
parameters $\mu_{2k-1}$ and $\mu_{2k+1}$ of different signs.
Obviously, if these parameters have the same sign, system (2.6)
has no limit cycles surrounding the origin at all.
    \par
Let system (2.6) have the unique limit cycle $\Gamma\!_{1}.$ Fix
the parameters $\mu_{2k+1}>0,$ $\mu_{2k-1}<0$ and input the third
parameter, $\mu_{2k-3}>0,$ into this system:
    $$
    \dot{x}=y,
    \quad
    \dot{y}=-x+y^{2}+\ldots+\mu_{2k-3}\,y^{2k-3}+y^{2k-2}
    +\ldots+\mu_{2k+1}\,y^{2k+1}.
    \eqno(2.8)
    $$
The vector field of~(2.8) is rotated counterclockwise, the focus
at the origin changes the character of its stability and the
second (unstable) limit cycle, $\Gamma_{2},$ immediately appears
from this point. Under further increasing $\mu_{2k-3},$ the limit
cycle $\Gamma_{2}$ will join with $\Gamma\!_{1}$ forming a
semi-stable limit cycle, $\Gamma_{\!12},$ which will disappear in
a ``trajectory concentration'' surrounding the origin. Can another
semi-stable limit cycle appear around the origin in addition to
$\Gamma_{\!12}?$ It is clear that such a limit cycle cannot appear
neither in the domain $D_{1}$ bounded on the inside by the cycle
$\Gamma\!_{1}$ nor in the domain $D_{3}$ bounded by the origin and
$\Gamma_{2}$ because of increasing the distance between the spiral
coils filling these domains under increasing the parameter
$\mu_{2k-3}.$
    \par
To prove impossibility of the appearance of a semi-stable limit
cycle in the domain $D_{2}$ bounded by the cycles $\Gamma\!_{1}$
and $\Gamma_{2}$ (before their joining), suppose the contrary,
i.\,e., for some set of values of the parameters,
$\mu_{2k+1}^{*}>0,$ $\mu_{2k-1}^{*}<0,$ and $\mu_{2k-3}^{*}>0,$
such a semi-stable cycle exists. Return to system (2.4) again and
input first the parameters $\mu_{2k-3}>0$ and $\mu_{2k+1}>0\!:$
    $$
    \dot{x}=y,
    \quad
    \dot{y}=-x+y^{2}+\ldots+\mu_{2k-3}\,y^{2k-3}
    +y^{2k-2}+y^{2k}+\mu_{2k+1}\,y^{2k+1}.
    \eqno(2.9)
    $$
Both parameters act in a similar way: they rotate the vector field
of (2.9) counterclockwise turning the origin into a nonrough
unstable focus.
    \par
Fix these parameters under $\mu_{2k-3}=\mu_{2k-3}^{*},$
$\mu_{2k+1}=\mu_{2k+1}^{*}$ and input the parameter $\mu_{2k-1}<0$
into (2.9) getting again system (2.8). Since, on our assumption,
this system has two limit cycles for $\mu_{2k-1}>\mu_{2k-1}^{*},$
there exists some value of the parameter, $\mu_{2k-1}^{12}$
$(\mu_{2k-1}^{*}<\mu_{2k-1}^{12}<0),$ for which a semi-stable
limit cycle, $\Gamma\!_{12},$ appears in system (2.8) and then
splits into a stable cycle, $\Gamma\!_{1},$ and an unstable cycle,
$\Gamma_{2},$ under further decreasing $\mu_{2k-1}.$ The formed
domain $D_{2}$ bounded by the limit cycles $\Gamma\!_{1},$
$\Gamma_{2}$ and filled by the spirals will enlarge since, on the
properties of a field rotation parameter, the interior unstable
limit cycle $\Gamma_{2}$ will contract and the exterior stable
limit cycle $\Gamma\!_{1}$ will expand under decreasing
$\mu_{2k-1}.$ The distance between the spirals of the domain
$D_{2}$ will naturally increase, what will prohibit from the
appearance of a semi-stable limit cycle in this domain for
$\mu_{2k-1}<\mu_{2k-1}^{12}.$
    \par
Thus, there are no such values of the parameters,
$\mu_{2k+1}^{*}>0,$ $\mu_{2k-1}^{*}<0,$ and $\mu_{2k-3}^{*}>0,$
for which system (2.8) would have an additional semi-stable limit
cycle. Obviously, there are no other values of the parameters
$\mu_{2k+1},$ $\mu_{2k-1},$ and $\mu_{2k-3}$ for which system
(2.8) would have more than two limit cycles surrounding the
origin. Therefore, two is the maximum number of limit cycles for
system (2.8). This result agrees with~\cite{rych}, where it was
proved for the first time that the maximum number of limit cycles
for Li\'{e}nard's system of the form
    $$
    \dot{x}=y,
    \quad
    \dot{y}=-x+\mu_{1}\,y+\mu_{3}\,y^{3}+\mu_{5}\,y^{5}
    \eqno(2.10)
    $$
was equal to two.
    \par
Suppose that system (2.8) has two limit cycles, $\Gamma\!_{1}$ and
$\Gamma_{2}$ (this is always possible if
$\mu_{2k+1}\gg-\mu_{2k-1}\gg\mu_{2k-3}>0),$ fix the parameters
$\mu_{2k+1},$ $\mu_{2k-1},$ $\mu_{2k-3}$ and consider a more
general system than (2.8) (and (2.10)) inputting the fourth
parameter, $\mu_{2k-5}<0,$ into (2.8):
    $$
    \dot{x}=y,
    \quad
    \dot{y}=-x+y^{2}+\ldots+\mu_{2k-5}\,y^{2k-5}
    +y^{2k-4}+\ldots+\mu_{2k+1}\,y^{2k+1}.
    \eqno(2.11)
    $$
Under decreasing $\mu_{2k-5},$ the vector field of~(2.11) will be
rotated clockwise and the focus at the origin will immediately
change the character of its stability generating the third
(stable) limit cycle, $\Gamma_{3}.$ Under further decreasing
$\mu_{2k-5},$ $\Gamma_{3}$ will join with $\Gamma_{2}$ forming a
semi-stable limit cycle, $\Gamma_{\!23},$ which will disappear in
a ``trajectory concentration'' surrounding the origin; the cycle
$\Gamma\!_{1}$ will expand infinitely tending to the Poincar\'{e}
circle at infinity.
    \par
Let system (2.11) have three limit cycles: $\Gamma\!_{1},$
$\Gamma_{2},$ $\Gamma_{3}.$ Could an additional semi-stable limit
cycle appear under decreasing $\mu_{2k-5},$ after splitting of
which system (2.11) would have five limit cycles around the
origin? It is clear that such a limit cycle cannot appear neither
in the domain $D_{2}$ bounded by the cycles $\Gamma\!_{1}$ and
$\Gamma_{2}$ nor in the domain $D_{4}$ bounded by the origin and
$\Gamma_{3}$ because of increasing the distance between the spiral
coils filling these domains under decreasing $\mu_{2k-5}.$
Consider two other domains: $D_{1}$ bounded on the inside by the
cycle $\Gamma\!_{1}$ and $D_{3}$ bounded by the cycles
$\Gamma_{2}$ and $\Gamma_{3}.$ As before, we will prove
impossibility of the appearance of a semi-stable limit cycle in
these domains by contradiction.
    \par
Suppose that for some set of values of the parameters
$\mu_{2k+1}^{*}>0,$ $\mu_{2k-1}^{*}<0,$ $\mu_{2k-3}^{*}>0,$ and
$\mu_{2k-5}^{*}<0,$ such a semi-stable cycle exists. Return to
system (2.4) again, input first the parameters $\mu_{2k-5}<0,$
$\mu_{2k-1}<0$ and then the parameter $\mu_{2k+1}>0\!:$
    $$
    \dot{x}=y,
    \quad
    \dot{y}=-x\!+\!y^{2}\!+\!\ldots\!+\!\mu_{2k-5}y^{2k-5}\!+\!\ldots\!
    +\!\mu_{2k-1}y^{2k-1}\!+\!y^{2k}\!+\!\mu_{2k+1}y^{2k+1}.
    \eqno(2.12)
    $$
Fix the parameters $\mu_{2k-5},$ $\mu_{2k-1}$ under the values
$\mu_{2k-5}^{*},$ $\mu_{2k-1}^{*},$ respectively. Under increasing
$\mu_{2k+1},$ the node at infinity will change the character of
its stability, the separatrix behaviour of the infinite saddle
will be also changed and a stable limit cycle, $\Gamma\!_{1},$
will immediately appear from the Poincar\'{e} circle at
infinity~\cite{BL}. Fix $\mu_{2k+1}$ under the value
$\mu_{2k+1}^{*}$ and input the parameter $\mu_{2k-3}>0$ into
(2.12) getting system (2.11).
    \par
Since, on our assumption, (2.11) has three limit cycles for
$\mu_{2k-3}<\mu_{2k-3}^{*},$ there exists some value of the
parameter $\mu_{2k-3}^{23}$ $(0<\mu_{2k-3}^{23}<\mu_{2k-3}^{*})$
for which a semi-stable limit cycle, $\Gamma_{23},$ appears in
this system and then splits into an unstable cycle, $\Gamma_{2},$
and a stable cycle, $\Gamma_{3},$ under further increasing
$\mu_{2k-3}.$ The formed domain $D_{3}$ bounded by the limit
cycles $\Gamma_{2},$ $\Gamma_{3}$ and also the domain $D_{1}$
bounded on the inside by the limit cycle $\Gamma\!_{1}$ will
enlarge and the spirals filling these domains will untwist
excluding a possibility of the appearance of a semi-stable limit
cycle there.
    \par
All other combinations of the parameters $\mu_{2k+1},$
$\mu_{2k-1},$ $\mu_{2k-3},$ and $\mu_{2k-5}$ are considered in a
similar way. It follows that system (2.11) has at most three limit
cycles. If we continue the procedure of successive inputting the
odd parameters, $\mu_{2k-7},\ldots,$ $\mu_{3},$ $\mu_{1},$ into
system (2.4), it is possible first to obtain $k$~limit cycles
$(\mu_{2k+1}\gg-\mu_{2k-1}\gg\mu_{2k-3}\gg-\mu_{2k-5}\gg\mu_{2k-7}\gg\ldots)$
and then to conclude that canonical system (2.1) (i.\,e.,
Li\'{e}nard's polynomial system (1.2) as well) has at most
$k$~limit cycles. The theorem is proved.

\section{An arbitrary polynomial system}
\label{3}

Let us consider an arbitrary polynomial system
    $$
    \dot{x}=P_{n}(x,y,\mu_{1},\ldots,\mu_{k}),
    \quad
    \dot{y}=Q_{n}(x,y,\mu_{1},\ldots,\mu_{k})
    \eqno(3.1)
    $$
containing $k$ field rotation parameters,
$\mu_{1},\ldots,\mu_{k},$ and having an anti-saddle at the origin.
Generalizing the main result of the previous section on the
maximum number of limit cycles surrounding a singular point in
Li\'{e}nard's polynomial system (1.2), we prove the following
theorem.
    \par
    \medskip
\noindent\textbf{Theorem 3.1.}
    \emph{Polynomial system $(3.1)$ containing $k$ field rotation
parameters and having a singular point of the center type at the
origin for the zero values of these parameters can have at most
$k-1$~limit cycles surrounding the origin.}
    \medskip
    \par
    \noindent\textbf{Proof.} \
Vanish all parameters of (3.1) and suppose that the obtained
system
    $$
    \dot{x}=P_{n}(x,y,0,\ldots,0),
    \quad
    \dot{y}=Q_{n}(x,y,0,\ldots,0)
    \eqno(3.2)
    $$
has a singular point of the center type at the origin. Let us
input successively the field rotation parameters,
$\mu_{1},\ldots,\mu_{k},$ into this system.
    \par
Suppose, for example, that $\mu_{1}>0$ and that the vector field
of the system
    $$
    \dot{x}=P_{n}(x,y,\mu_{1},0,\ldots,0),
    \quad
    \dot{y}=Q_{n}(x,y,\mu_{1},0,\ldots,0)
    \eqno(3.3)
    $$
is rotated counterclockwise turning the origin into a stable focus
under increasing $\mu_{1}.$
    \par
Fix $\mu_{1}$ and input the parameter $\mu_{2}$ into (3.3)
changing it so that the field of the system
    $$
    \dot{x}=P_{n}(x,y,\mu_{1},\mu_{2},0,\ldots,0),
    \quad
    \dot{y}=Q_{n}(x,y,\mu_{1},\mu_{2},0,\ldots,0)
    \eqno(3.4)
    $$
would be rotated in opposite direction (clockwise). Let be so for
$\mu_{2}<0.$ Then, for some value of this parameter, a limit cycle
will appear in system (3.4). There are three logical possibilities
for such a bifurcation: 1)~the limit cycle appears from the focus
at the origin; 2)~it can also appear from some separatrix cycle
surrounding the origin; 3)~the limit cycle appears from a
so-called ``trajectory concentration''. In the last case, the
limit cycle is semi-stable and, under further decreasing
$\mu_{2},$ it splits into two limit cycles (stable and unstable),
one of which then disappears at (or tends to) the origin and the
other disappears on (or tends to) some separatrix cycle
surrounding this point. But since the stability character of both
a singular point and a separatrix cycle is quite easily
controlled~\cite{Gaiko}, this logical possibility can be excluded.
Let us choose one of the two other possibilities: for example, the
first one, the so-called Andronov--Hopf bifurcation. Suppose that,
for some value of~$\mu_{2},$ the focus at the origin becomes
non-rough, changes the character of its stability and generates a
stable limit cycle, $\Gamma\!_{1}.$
    \par
Under further decreasing $\mu_{2},$ three new logical
possibilities can arise: 1)~the limit cycle $\Gamma\!_{1}$
disappears on some separatrix cycle surrounding the origin;
2)~a~sepa\-ratrix cycle can be formed earlier than $\Gamma\!_{1}$
disappears on it, then it generates one more (unstable) limit
cycle, $\Gamma_{2},$ which joins with $\Gamma\!_{1}$ forming a
semi-stable limit cycle, $\Gamma\!_{12},$ disappearing in a
``trajectory concentration'' under further decreasing $\mu_{2};$
3)~in the domain $D_{1}$ outside the cycle $\Gamma\!_{1}$ or in
the domain $D_{2}$ inside $\Gamma\!_{1},$ a semi-stable limit
cycle appears from a ``trajectory concentration'' and then splits
into two limit cycles (logically, the appearance of such
semi-stable limit cycles can be repeated).
    \par
Let us consider the third case. It is clear that, under decreasing
$\mu_{2},$ a semi-stable limit cycle cannot appear in the domain
$D_{2},$ since the focus spirals filling this domain will untwist
and the distance between their coils will increase because of the
vector field rotation. By contradiction, we can prove that a
semi-stable limit cycle cannot appear in the domain $D_{1}.$
Suppose it appears in this domain for some values of the
parameters $\mu_{1}^{*}>0$ and $\mu_{2}^{*}<0.$ Return to initial
system (3.2) and change the inputting order for the field rotation
parameters. Input first the parameter $\mu_{2}<0\!:$
    $$
    \dot{x}=P_{n}(x,y,\mu_{2},0,\ldots,0),
    \quad
    \dot{y}=Q_{n}(x,y,\mu_{2},0,\ldots,0).
    \eqno(3.5)
    $$
Fix it under $\mu_{2}=\mu_{2}^{*}.$ The vector field of~(3.5) is
rotated clockwise and the origin turns into a unstable focus.
Inputting the parameter $\mu_{1}>0$ into (3.5), we get again
system (3.4), the vector field of which is rotated
counterclockwise. Under this rotation, a stable limit cycle,
$\Gamma\!_{1},$ will appear from some separatrix cycle. The limit
cycle $\Gamma\!_{1}$ will contract, the outside spirals winding
onto this cycle will untwist and the distance between their coils
will increase under increasing $\mu_{1}$ to the value
$\mu_{1}^{*}.$ It follows that there are no values of
$\mu_{2}^{*}<0$ and $\mu_{1}^{*}>0,$ for which a semi-stable limit
cycle could appear in the domain $D_{1}.$
    \par
The second logical possibility can be excluded by controlling the
stability character of the separatrix cycle~\cite{Gaiko}. Thus,
only the first possibility is valid, i.\,e., system (3.4) has at
most one limit cycle.
    \par
Let system (3.4) have the unique limit cycle $\Gamma\!_{1}.$ Fix
the parameters $\mu_{1}>0,$ $\mu_{2}<0$ and input the third
parameter, $\mu_{3}>0,$ into this system supposing that $\mu_{3}$
rotates its vector field counterclockwise:
    $$
    \dot{x}=P_{n}(x,y,\mu_{1},\mu_{2},\mu_{3},0,\ldots,0),
    \quad
    \dot{y}=Q_{n}(x,y,\mu_{1},\mu_{2},\mu_{3},0,\ldots,0).
    \eqno(3.6)
    $$
Here we can have two basic possibilities: 1)~the limit cycle
$\Gamma\!_{1}$ disappears at the origin; 2)~the second (unstable)
limit cycle, $\Gamma_{2},$ appears from the origin and, under
further increasing the parameter $\mu_{3},$ the cycle $\Gamma_{2}$
joins with $\Gamma\!_{1}$ forming a semi-stable limit cycle,
$\Gamma_{\!12},$ which disappears in a ``trajectory
concentration'' surrounding the origin. Besides, we can also
suggest that: 3)~in the domain $D_{2}$ bounded by the origin and
$\Gamma\!_{1},$ a semi-stable limit cycle, $\Gamma_{23},$ appears
from a ``trajectory concentration'', splits into an unstable
cycle, $\Gamma_{2},$ and a stable cycle, $\Gamma_{3},$ and then
the cycles $\Gamma\!_{1},$ $\Gamma_{2}$ disappear through a
semi-stable limit cycle, $\Gamma\!_{12},$ and the cycle
$\Gamma_{3}$ disappears through the Andronov--Hopf bifurcation;
4)~a~semi-stable limit cycle, $\Gamma_{34},$ appears in the domain
$D_{2}$ bounded by the cycles $\Gamma\!_{1},$ $\Gamma_{2}$ and,
for some set of values of the parameters, $\mu_{1}^{*},$
$\mu_{2}^{*},$ $\mu_{3}^{*},$ system (3.6) has at least four limit
cycles.
    \par
Let us consider the last, fourth, case. It is clear that a
semi-stable limit cycle cannot appear neither in the domain
$D_{1}$ bounded on the inside by the cycle $\Gamma\!_{1}$ nor in
the domain $D_{3}$ bounded by the origin and $\Gamma_{2}$ because
of increasing the distance between the spiral coils filling these
domains under increasing the parameter $\mu_{3}.$ To prove
impossibility of the appearance of a semi-stable limit cycle in
the domain $D_{2},$ suppose the contrary, i.\,e., for some set of
values of the parameters, $\mu_{1}^{*}>0,$ $\mu_{2}^{*}<0,$ and
$\mu_{3}^{*}>0,$ such a semi-stable cycle exists. Return to system
(3.2) again and input first the parameters $\mu_{3}>0,$
$\mu_{1}>0\!:$
    $$
    \dot{x}=P_{n}(x,y,\mu_{1},\mu_{3},0,\ldots,0),
    \quad
    \dot{y}=Q_{n}(x,y,\mu_{1},\mu_{3},0,\ldots,0).
    \eqno(3.7)
    $$
Fix these parameters under $\mu_{3}=\mu_{3}^{*},$
$\mu_{1}=\mu_{1}^{*}$ and input the parameter $\mu_{2}<0$ into
(3.7) getting again system (3.6). Since, on our assumption, this
system has two limit cycles for $\mu_{2}>\mu_{2}^{*},$ there
exists some value of the parameter, $\mu_{2}^{12}$
$(\mu_{2}^{*}<\mu_{2}^{12}<0),$ for which a semi-stable limit
cycle, $\Gamma\!_{12},$ appears in system (3.6) and then splits
into a stable cycle, $\Gamma\!_{1},$ and an unstable cycle,
$\Gamma_{2},$ under further decreasing $\mu_{2}.$ The formed
domain $D_{2}$ bounded by the limit cycles $\Gamma\!_{1},$
$\Gamma_{2}$ and filled by the spirals will enlarge, since, on the
properties of a field rotation parameter, the interior unstable
limit cycle $\Gamma_{2}$ will contract and the exterior stable
limit cycle $\Gamma\!_{1}$ will expand under decreasing $\mu_{2}.$
The distance between the spirals of the domain $D_{2}$ will
naturally increase, what will prohibit from the appearance of a
semi-stable limit cycle in this domain for $\mu_{2}<\mu_{2}^{12}.$
    \par
Thus, there are no such values of the parameters, $\mu_{1}^{*}>0,$
$\mu_{2}^{*}<0,$ $\mu_{3}^{*}>0,$ for which system (3.6) would
have an additional semi-stable limit cycle. Therefore, the fourth
case cannot be realized. The third case is considered absolutely
similarly. It follows from the first two cases that system (3.6)
can have at most two limit cycles.
    \par
Suppose that system (3.6) has two limit cycles, $\Gamma\!_{1}$ and
$\Gamma_{2},$ fix the parameters $\mu_{1}>0,$ $\mu_{2}<0,$
$\mu_{3}>0$ and input the fourth parameter, $\mu_{4}<0,$ into this
system supposing that $\mu_{4}$ rotates its vector field
clockwise:
    $$
    \dot{x}=P_{n}(x,y,\mu_{1},\ldots,\mu_{4},0,\ldots,0),
    \quad
    \dot{y}=Q_{n}(x,y,\mu_{1},\ldots,\mu_{4},0,\ldots,0).
    \eqno(3.8)
    $$
The most interesting logical possibility here is that when the
third (stable) limit cycle, $\Gamma_{3},$ appears from the origin
and then, under preservation of the cycles $\Gamma\!_{1}$ and
$\Gamma_{2},$ in the domain $D_{3}$ bounded on the inside by the
cycle $\Gamma_{3}$ and on the outside by the cycle $\Gamma_{2},$ a
semi-stable limit cycle, $\Gamma_{\!45},$ appears and then splits
into a stable cycle, $\Gamma\!_{4},$ and an unstable cycle,
$\Gamma_{5},$ i.\,e., when system (3.8) for some set of values of
the parameters, $\mu_{1}^{*},$ $\mu_{2}^{*},$ $\mu_{3}^{*},$
$\mu_{4}^{*},$ has at least five limit cycles. Logically, such a
semi-stable limit cycle could also appear in the domain $D_{1}$
bounded on the inside by the cycle $\Gamma\!_{1},$ since, under
decreasing $\mu_{4},$ the spirals of the trajectories of (3.8)
will twist and the distance between their coils will decrease. On
the other hand, in the domain $D_{2}$ bounded on the inside by the
cycle $\Gamma_{2}$ and on the outside by the cycle $\Gamma\!_{1}$
and also in the domain $D_{4}$ bounded by the origin and
$\Gamma_{3},$ a semi-stable limit cycle cannot appear, since,
under decreasing $\mu_{4},$ the spirals will untwist and the
distance between their coils will increase. To prove impossibility
of the appearance of a semi-stable limit cycle in the domains
$D_{3}$ and $D_{1},$ suppose the contrary, i.\,e., for some set of
values of the parameters, $\mu_{1}^{*}>0,$ $\mu_{2}^{*}<0,$
$\mu_{3}^{*}>0,$ and $\mu_{4}^{*}<0,$ such a semi-stable cycle
exists. Return to system (3.2) again, input first the parameters
$\mu_{4}<0,$ $\mu_{2}<0$ and then the parameter $\mu_{1}>0\!:$
    $$
    \dot{x}=P_{n}(x,y,\mu_{1},\mu_{2},\mu_{4},0,\ldots,0),
    \quad
    \dot{y}=Q_{n}(x,y,\mu_{1},\mu_{2},\mu_{4},0,\ldots,0).
    \eqno(3.9)
    $$
Fix the parameters $\mu_{4},$ $\mu_{2}$ under the values
$\mu_{4}^{*},$ $\mu_{2}^{*},$ respectively. Under increasing
$\mu_{1},$ a separatrix cycle is formed around the origin
generating a stable limit cycle, $\Gamma\!_{1}.$ Fix $\mu_{1}$
under the value $\mu_{1}^{*}$ and input the parameter $\mu_{3}>0$
into (3.9) getting system (3.8).
    \par
Since, on our assumption, system (3.8) has three limit cycles for
$\mu_{3}<\mu_{3}^{*},$ there exists some value of the parameter
$\mu_{3}^{23}$ $(0<\mu_{3}^{23}<\mu_{3}^{*})$ for which a
semi-stable limit cycle, $\Gamma_{23},$ appears in this system and
then splits into an unstable cycle, $\Gamma_{2},$ and a stable
cycle, $\Gamma_{3},$ under further increasing $\mu_{3}.$ The
formed domain $D_{3}$ bounded by the limit cycles $\Gamma_{2},$
$\Gamma_{3}$ and also the domain $D_{1}$ bounded on the inside by
the limit cycle $\Gamma\!_{1}$ will enlarge and the spirals
filling these domains will untwist excluding a possibility of the
appearance of a semi-stable limit cycle there.
    \par
All other combinations of the parameters $\mu_{1},$ $\mu_{2},$
$\mu_{3},$ and $\mu_{4}$ are considered in a similar way. It
follows that system (3.8) has at most three limit cycles. If we
continue the procedure of successive inputting the field rotation
parameters, $\mu_{5},$ $\mu_{6},\ldots,$ $\mu_{k},$ into system
(3.2), it is possible to conclude that system (3.1) can have at
most $k-1$~limit cycles surrounding the origin. The theorem is
proved.

\section{A generalized Li\'{e}nard's system}
\label{4}

In \cite{gaivh}, we considered a generalized Li\'{e}nard's cubic
system of the form:
    $$
    \dot{x}=y,\quad
    \displaystyle\dot{y}=-x+(\lambda-\mu)\,y+(3/2)\,x^2
        +\mu\,xy-(1/2)\,x^3+\alpha\,x^2y\,.
    \eqno(4.1)
    $$
This system has three finite singularities: a saddle $(1,0)$ and
two antisaddles~--- $(0,0)$ and $(2,0).$ At infinity system
$(4.1)$ can have either the only nilpotent singular point of
fourth order with two closed elliptic and four hyperbolic domains
or two singular points: one of them is a hyperbolic saddle and the
other is a triple nilpotent singular point with two elliptic and
two hyperbolic domains. We studied global bifurcations of limit
and separatrix cycles of (4.1), found possible distributions of
its limit cycles and carried out a classification of its
separatrix cycles. We proved also the following theorems.
    \par
    \medskip
\noindent \textbf{Theorem 4.1.} \
    \emph{The foci of system $(4.1)$ can be at most of second order.}
    \par
\noindent \textbf{Theorem 4.2.} \
    \emph{System~$(4.1)$ has at least three limit cycles.}
    \medskip
    \par
Using the results obtained in~\cite{gaivh} and applying the
approach developed in this paper, we can easily prove a much
stronger theorem.
    \par
    \medskip
\noindent \textbf{Theorem 4.3.} \
    \emph{System~$(4.1)$ has at most three limit cycles with the
following their distributions$\,:$ $((1,1),1),$ $((1,2),0),$
$((2,1),0),$ $((1,0),2),$ $((0,1),2),$ where the first two numbers
denote the numbers of limit cycles surrounding each of two
anti-saddles and the third one denotes the number of limit cycles
surrounding simultaneously all three finite singularities.}
    \par
    \medskip
Theorem 4.3 agrees, for example, with the earlier results by Iliev
and\linebreak Perko~\cite{ilievp}, but it does not agree with a
quite recent result by Dumortier and Li~\cite{dumli} published in
the same journal. The authors of both papers use very similar
methods: small perturbations of a Hamiltonian system.
In~\cite{ilievp}, the zeros of the Melnikov functions are studied
and, in particular, it is proved that at most two limit cycles can
bifurcate from either the interior or exterior period annulus of
the Hamiltonian under small parameter perturbations giving a
ge\-ne\-ra\-lized Li\'{e}nard system. In~\cite{dumli}, zeros of
the Abelian integrals are studied and it is ``proved'' that at
most four limit cycles can bifurcate from the exterior period
annulus. Thus, Dumortier and Li ``obtain'' a configuration of four
big limit cycles surrounding three finite singularities together
with the fifth small limit cycle which surrounds one of the
anti-saddles.
    \par
The result by Dumortier and Li~\cite{dumli} also does not agree
with the Wintner--Perko termination principle for multiple limit
cycles~\cite{Gaiko}, \cite{Perko}. Applying the method as
developed in~\cite{bgai}, \cite{gai1}--\cite{gaib}, we can show
that system~(4.1) cannot have neither a multiplicity-three limit
cycle nor more than three limit cycles in any configuration. That
will be another proof of Theorem~4.3 (the same approach can be
applied to proving Theorems~2.2 and~3.1 as well). But first let us
formulate the Wintner--Perko termination principle~\cite{Perko}
for the polynomial system
    $$
    \mbox{\boldmath$\dot{x}$}=\mbox{\boldmath$f$}
    (\mbox{\boldmath$x$},\mbox{\boldmath$\mu$)},
    \eqno(4.2_{\mbox{\boldmath$\mu$}})
    $$
where $\mbox{\boldmath$x$}\in\textbf{R}^2;$ \
$\mbox{\boldmath$\mu$}\in\textbf{R}^n;$ \
$\mbox{\boldmath$f$}\in\textbf{R}^2$ \ $(\mbox{\boldmath$f$}$ is a
polynomial vector function).
    \par
    \medskip
\noindent\textbf{Theorem 4.4
    (Wintner--Perko termination principle).}
    \emph{Any one-para\-me\-ter fa\-mi\-ly of multiplicity-$m$
limit cycles of relatively prime polynomial system \
$(4.2_{\mbox{\boldmath$\mu$}})$ can be extended in a unique way to
a maximal one-parameter family of multiplicity-$m$ limit cycles of
\ $(4.2_{\mbox{\boldmath$\mu$}})$ which is either open or cyclic.}
    \par
\emph{If it is open, then it terminates either as the parameter or
the limit cycles become unbounded; or, the family terminates
either at a singular point of \ $(4.2_{\mbox{\boldmath$\mu$}}),$
which is typically a fine focus of multiplicity~$m,$ or on a
$($compound$\,)$ separatrix cycle of \
$(4.2_{\mbox{\boldmath$\mu$}}),$ which is also typically of
multiplicity~$m.$}
    \medskip
    \par
The proof of this principle for general polynomial system
$(4.2_{\mbox{\boldmath$\mu$}})$ with a vector parameter
$\mbox{\boldmath$\mu$}\in\textbf{R}^n$ parallels the proof of the
pla\-nar termination principle for the system
    $$
    \vspace{1mm}
    \dot{x}=P(x,y,\lambda),
        \quad
    \dot{y}=Q(x,y,\lambda)
    \eqno(4.2_{\lambda})
    \vspace{2mm}
    $$
with a single parameter $\lambda\in\textbf{R}$ (see \cite{Gaiko},
\cite{Perko}), since there is no loss of generality in assuming
that system $(4.2_{\mbox{\boldmath$\mu$}})$ is parameterized by a
single parameter $\lambda;$ i.\,e., we can assume that there
exists an analytic mapping $\mbox{\boldmath$\mu$}(\lambda)$ of
$\textbf{R}$ into $\textbf{R}^n$ such that
$(4.2_{\mbox{\boldmath$\mu$}})$ can be written as
$(4.2\,_{\mbox{\boldmath$\mu$}(\lambda)})$ or even
$(4.2_{\lambda})$ and then we can repeat everything, what had been
done for system $(4.2_{\lambda})$ in~\cite{Perko}. In particular,
if $\lambda$ is a field rotation parameter of $(4.2_{\lambda}),$
it is valid the following Perko's theorem on monotonic families of
limit cycles.
    \par
    \medskip
\noindent\textbf{Theorem 4.5.}
    \emph{If $L_{0}$ is a nonsingular multiple limit cycle of
$(4.2_{0}),$ then  $L_{0}$ belongs to a one-parameter family of
limit cycles of $(4.2_{\lambda});$ furthermore\/$:$}
    \par
1)~\emph{if the multiplicity of $L_{0}$ is odd, then the family
either expands or contracts mo\-no\-to\-ni\-cal\-ly as $\lambda$
increases through $\lambda_{0};$}
    \par
2)~\emph{if the multiplicity of $L_{0}$ is even, then $L_{0}$
be\-fur\-cates into a stable and an unstable limit cycle as
$\lambda$ varies from $\lambda_{0}$ in one sense and $L_{0}$
dis\-ap\-pears as $\lambda$ varies from $\lambda_{0}$ in the
opposite sense; i.\,e., there is a fold bifurcation at
$\lambda_{0}.$}
    \medskip
    \par
\noindent\textbf{Proof of Theorem 4.3.} \ The proof is carried out
by contradiction. Suppose that system (4.1) with three field
rotation parameters, $\lambda,$ $\mu,$ and $\alpha,$ has three
limit cycles around, for example, the origin (the case when limit
cycles surround another focus is considered in a similar way).
Then we get into some domain in the space of these parameters
which is bounded by two fold bifurcation surfaces forming a cusp
bifurcation surface of multiplicity-three limit cycles.
    \par
The cor\-res\-pon\-ding maximal one-parameter family of
multiplicity-three limit cycles cannot be cyclic, otherwise there
will be at least one point cor\-res\-pon\-ding to the limit cycle
of multi\-pli\-ci\-ty four (or even higher) in the parameter
space. Extending the bifurcation curve of multi\-pli\-ci\-ty-four
limit cycles through this point and parameterizing the
corresponding maximal one-parameter family of
multi\-pli\-ci\-ty-four limit cycles by a field-rotation
para\-me\-ter, according to Theorem~4.5, we will obtain a
monotonic curve which, by the Wintner--Perko termination principle
(Theorem~4.4), terminates either at the origin or on some
separatrix cycle surrounding the origin. Since we know absolutely
precisely at least the cyclicity of the singular point
(Theorem~4.1) which is equal to two, we have got a contradiction
with the termination principle stating that the multiplicity of
limit cycles cannot be higher than the multi\-pli\-ci\-ty
(cyclicity) of the singular point in which they terminate.
    \par
If the maximal one-parameter family of multiplicity-three limit
cycles is not cyclic, on the same principle (Theorem~4.4), this
again contradicts to Theorem~4.1 not admitting the multiplicity of
limit cycles higher than two. Moreover, it also follows from the
termination principle that neither the ordinary separatrix loop
nor the eight-loop cannot have the multiplicity (cyclicity) higher
than two (in that way, it can be proved that the cyclicity of
three other separatrix cycles~\cite{gaivh} is at most two).
Therefore, according to the same principle, there are no more than
two limit cycles in the exterior domain surrounding all three
finite singularities of~(4.1). Thus, system~(4.1) cannot have
neither a multiplicity-three limit cycle nor more than three limit
cycles in any configuration. The theorem is proved.
    \par
So, we have found two approaches to solving \emph{Smale's
Thirteenth and Hil\-bert's Sixteenth Problems.} Both these
approaches are based on the application of field rotation
parameters which determine limit cycle bifurcations of polynomial
systems.

\end{document}